\documentclass[12pt]{amsart}
\usepackage{mathptmx}
\usepackage{amsmath}
\usepackage{amssymb}
\usepackage{array}
\usepackage{geometry}
\usepackage[bookmarks=true,colorlinks=true, pdfstartview=FitV, linkcolor=black, citecolor=blue, urlcolor=black]{hyperref}
\usepackage{movie15}
\usepackage{graphicx}
\usepackage{epsfig}
\usepackage{amssymb,amsthm,amsmath,amsfonts,amstext,mathrsfs}
\usepackage{latexsym}
\usepackage{url}
\input xy
\xyoption{all}

\begin{document}

\title{\bf {On the affirmative solution to Salem's problem}}
\author{Semyon Yakubovich}
\maketitle

 \markboth{\rm \centerline{SEMYON YAKUBOVICH}}{}
 \markright{\rm \centerline{SALEM-TYPE PROBLEM}}

 \begin{abstract}{The  Salem problem to verify whether Fourier-Stieltjes coefficients of the Minkowski question mark function vanish at infinity is solved recently affirmatively. In this paper by using methods of classical analysis and special functions  we solve a Salem-type problem about the behavior at infinity of a linear combination of the Fourier-Stieltjes transforms.   Moreover, as a consequence of the Salem problem, some  asymptotic relations at infinity for  the Fourier-Stieltjes coefficients of a  power $m\in \mathbb{N}$  of the Minkowski question mark function are derived.}
\end{abstract}



\section{Introduction}

Let $x \in \mathbb{R}$ and consider  the following Fourier-Stieltjes transforms
$$f (x)= \int_0^1  e^{ixt} d q(t),\eqno(1.1)$$
$$F(x)= \int_0^\infty e^{ixt} d q(t).\eqno(1.2)$$
Here $q(x)$ is the famous Minkowski  question mark function $?(x) \equiv q(x)$.   This  function is defined by \cite{D} $q(x):[0,1]\mapsto[0,1]$
$$q ([0,a_{1},a_{2},a_{3},\ldots])=2\sum\limits_{i=1}^{\infty}(-1)^{i+1}2^{-\sum_{j=1}^{i}a_{j}},$$
where $x=[0,a_{1},a_{2},a_{3},\ldots]$ stands for the representation of $x$ by a regular continued fraction.  It is well known that $q(x)$ is continuous, strictly increasing and supports a singular  measure. It satisfies  the  following functional  equations (cf. \cite{Al}, p.3), which will be used in the sequel
$$q (x)= 1- q(1-x),\  x \in [0,1],\eqno(1.3)$$
$$q(x)= 2 q\left(\frac{x}{x+1}\right),\  x >0,\eqno(1.4)$$
$$q(x)+ q\left({1\over x}\right)= 2, \  x > 0.\eqno(1.5)$$
When $x \to 0$, it decreases exponentially $ \  q(x)= O\left(2^{-
1/x}\right)$. Key values are $q(0)=0, \ q(1)= 1,\
q(\infty)= 2$. For instance, from (1.3)  and asymptotic
behavior of the Minkowski function near zero one can easily get the
finiteness of the following integrals
$$\int_0^1 x^\lambda\  d q(x) < \infty, \ \lambda \in \mathbb{R},\eqno(1.6)$$
$$\int_0^1 (1-x)^\lambda\ d q(x) < \infty, \ \lambda \in
\mathbb{R}.\eqno(1.7)$$
Further, as was proved by Salem \cite{Sal1}, the Minkowski question
mark function satisfies the H\"older condition
$$\left|q(x)-  q(y)\right| < C |x-y|^\alpha, \ \alpha < 1,$$
where
$$\alpha= \frac{\log 2}{2 \log {\sqrt 5 + 1\over 2}}= 0,7202_+.$$
and $C >0$ is an absolute constant.   As  we observe from the  functional equation (1.3)  the Fourier-Stieltjes transform (1.1)   satisfies the functional relation
$$f(x) = e^{ix} f(-x), \eqno(1.8)$$
and  therefore $e^{-ix/2} f(x)$ is real-valued. So, taking its imaginary part, we  obtain the equality
$$\cos \left({x\over 2}\right) f_{s} (x) = \sin \left({x\over 2}\right) f_c (x),\eqno(1.9)$$ 
where $f_{s},\ f_{c} $ are the Fourier-Stieltjes sine and cosine transforms of the Minkowski question mark function, respectively,
$$f_{s} (x)= \int_0^1  \sin(xt) d q(t),\eqno(1.10)$$
$$f_{c} (x)= \int_0^1  \cos(xt) d q(t).\eqno(1.11)$$
Hence, letting, for instance,  $x = 2\pi  n, \ n \in \mathbb{N}_0$ it
gives $f_{s}(2\pi n)=0$ and $f_{c} (2\pi n)= d_n$. In 1943  Salem asked \cite{Sal1}
whether $d_{n}\rightarrow 0$, as $n\rightarrow\infty$.   This question is quite delicate, since it concerns singular functions  (see \cite{Sal3}, Ch. IV)  and  the classical Riemann-Lebesgue lemma for the class $L_1$, in general, cannot be applied.  A singular function is defined as a continuous, bounded  monotone function with a null derivative almost everywhere. Hence it supports a positive bounded Borel measure, which is singular with respect to Lebesgue measure. For such singular measures there are various examples whose Fourier transforms do not tend to zero, although some do (see, for instance, in \cite{Sal1},  \cite{Sal2},  \cite{Men}). In \cite{Win} (see also \cite{IV}) it was proved that for every $\varepsilon >0$ there exists a singular monotone function, which supports a measure whose
Fourier-Stieltjes  transform behaves as $O(t^{-{1\over 2} +\varepsilon}), \ |t| \to \infty$.  

 In fact, it is worth to mention that the Salem problem  is an old and quite attractive problem in the number theory and Fourier analysis \cite{Walter}.  Several attempts were undertaken to solve Salem's problem by analytic methods (see, for instance,  in \cite{Al}, \cite{Yakeq}, \cite{YAK3} ).   Finally,  it was solved affirmatively in \cite{Jor} as  a special case of the general theory of Fourier transforms, involving  Gibbs measures for the Gauss map.

  In the sequel we will give  an affirmative  solution to a Salem-type  problem, using the methods of classical analysis and special functions.  It involves the asymptotic behavior at infinity of a linear combination of the Fourier-Stieltjes transforms.    Besides, in Section 2 we will establish  a new integro-differential equation for the Fourier-Stieltjes transform (1.1).  It involves, in turn,  the following functional equation, which is proved    by the author in \cite{Yakeq} (see Lemma 1 below) and relates to  transforms (1.1), (1.2)
$$f(x)= \left(1- {e^{ix}\over 2}\right) F(x),\ x \in \mathbb{R}.\eqno(1.12)$$   
Taking   real and imaginary parts of both sides in (1.12), we derive interesting equalities (see details in \cite{Yakeq}), which will be used below, namely
$$\int_{1}^{\infty}\cos xt \ d q(t)= \frac{1- 8\sin^2(x/2)}{1+ 8\sin^2(x/2)}\int_{0}^{1}\cos x t\ d
q(t),\  x \in \mathbb{R},\eqno(1.13)$$
$$\int_{1}^{\infty}\sin xt \ d q(t)= \frac{5- 8\sin^2(x/2)}{1+ 8\sin^2(x/2)}\int_{0}^{1}\sin  xt\ d
q(t),\ x \in \mathbb{R}.\eqno(1.14)$$
Making $x \to 0$ in (1.14), we find, in particular,  
$$\int_{1}^{\infty} t  d q(t)=  5\int_{0}^{1} t  dq(t).\eqno(1.15)$$
Moreover, using the functional equation (1.3), it can be proved the important equality  for coefficients $d_n$
 $$d_n= 2\int_{0}^{1} t \cos (2\pi n t) \ d q(t).\eqno(1.16)$$
Indeed, we have 
$$ \int_{0}^{1} t \cos (2\pi n t )\ d q(t) -  {d_n\over 2} = \left(\int_{0}^{1/2} + \int_{1/2}^{1}\right)  \left( t - {1\over 2} \right) \cos (2\pi n t)  \ d q(t) $$
$$=   \int_{0}^{1/2}  \left( t - {1\over 2} \right) \cos (2\pi n t)  \ d q(t) -  \int_{1/2}^{1} \left( t - {1\over 2} \right) \cos (2\pi n t)  \ d q(1-t) = 0.$$
Further, we give values of  the important integrals, which  will be employed in the sequel.  Precisely, according to relation (2.16.48.20) in \cite{Prud}, Vol. 2 and the differentiation with respect to a parameter,  the following integral with respect to an index of the modified Bessel function $K_{i\tau}(x)$ \cite{Olver} is calculated (cf. \cite{YAK3} )
$$ {1\over \pi}\int_{-\infty}^\infty \tau e^{\lambda \tau} \
\left(x+ (1+x^2)^{1/2}\right)^{i\tau} K_{i\tau}(t) d\tau\\
= t \exp\left(-t \left[(1+x^2)^{1/2}\cos \lambda - i x
\sin\lambda\right]\right)$$

$$\times  \left[(1+x^2)^{1/2}\sin \lambda + ix \cos
\lambda\right], \quad  x, t  > 0,  \quad  0  \le \lambda < {\pi \over 2}.\eqno(1.17)$$
Meanwhile, the Fourier cosine transform of the modified Bessel function   $K_{i\tau}(x)$ is given via relation (2.16.14.1)  in \cite{Prud}, Vol. 2

$$\int_0^\infty \cos( xy) \  K_{i\tau}(x)\ dx = {\pi\over 2} \  \frac{ \cos( \tau \log( y+ (y^2+1)^{1/2} )) } { (y^2+1)^{1/2} \ \cosh(\pi\tau/2)},\ y > 0,\ \tau \in \mathbb{R}.\eqno(1.18) $$
Finally, we will need the value of the integral  (cf. relation (2.4.4.3) in \cite{Prud}, Vol. 1)

$$\int_0^\infty {\cosh(a y)\over \cosh (\pi y/2)} dy = \sec( a), \quad |{\rm Re} a | < {\pi\over 2}.\eqno(1.19)$$
and the following generalization of the Dirichlet integral (cf. relations (2.5.25.13) in \cite{Prud}, Vol. 1)

$$\int_0^\infty y \sin ( b y  ) \  \cos( c (y^2+1)^{1/2} ) {dy \over y^2 + a^2} = {\pi\over 2}\  e^{-ab} \cos( c (1 -a^2)^{1/2} ),\quad  b >c >0, \ a >0. \eqno(1. 20)$$   

\section{Integro-differential equation for the Fourier-Stieltjes transform $(1.1)$}

In order to make the paper self-contained we  begin with the proof of the relation (1.12) (cf. \cite{Yakeq}).   

{\bf Lemma 1}.  {\it Let $x \in \mathbb{R}$ and $f(x)$, $F(x)$ be Fourier-Stieltjes transforms $(1.1),\ (1.2)$, respectively.  Then functional equation $(1.12)$ holds. }

\begin{proof} The proof is based on functional equations (1.4), (1.5) for the Minkowski question mark function and  simple properties of the Stieltjes integral.   In fact, we derive the chain of equalities
\begin{eqnarray*}
\int_{0}^{1}e^{ixt}\ d q(t)=
\int_{0}^{\infty}e^{ixt}\
d q(t)- \int_{1}^{\infty}e^{ixt}\ d q(t)\\
= \int_{0}^{\infty}e^{ixt}\ d q(t) - e^{ix}
\int_{0}^{\infty}e^{ixt}\
d q \left({ t+1}\right)\\
= \int_{0}^{\infty}e^{ixt}\ d q(t) + e^{ix}
\int_{0}^{\infty}e^{ixt}\
d q \left({1\over { t+1}}\right)\\
= \int_{0}^{\infty}e^{ixt}\ d q(t) + e^{ix}
\int_{0}^{\infty}e^{ixt}\
d  q\left({1/t \over { 1+ 1/t}}\right)\\
= \int_{0}^{\infty}e^{ixt}\ d q(t) + {e^{ix}\over 2}
\int_{0}^{\infty}e^{ixt}\
d q\left({1\over t}\right)\\
= \left(1- {e^{ix}\over 2}\right) \int_{0}^{\infty}e^{ixt}\ d
q(t),
\end{eqnarray*}
which yields (1.12).
\end{proof}

{\bf Theorem 1.}  {\it Let $ x \in \mathbb{R}_+$.  The Fourier-Stieltjes transform $(1.1)$ satisfies the following integro-differential equation, involving the operator of the modified Hankel transform}
$$\frac{ e^{ix}}{ 2 - e^{ix}  } \left[  f^\prime (x) + \frac{ 2 i   f(x)}{ 2 - e^{ix}  }\right] =  - \int_0^\infty J_0\left(2\sqrt {x y}\right)  e^{-iy}  f(y) dy.\eqno(2.1)$$

\begin{proof}  Indeed, differentiating (1.12) with respect to $x$ and using it again, we find 
$$f^\prime (x)= -   \frac{  i \  e^{ix}  f(x)}{ 2 - e^{ix}  } + i \left(1- {e^{ix}\over 2}\right) 
\int_0^\infty t e^{ixt} \  d q(t),\eqno(2.2)$$
where the differentiation under the integral sign in (1.2) is allowed via the  simple  estimate 
$$\left| \int_0^\infty t e^{ixt} \  d q(t)\right| \le   \int_0^\infty t  d q(t) =  3, $$ 
where the latter equality is due to (1.15), (1.16).  Hence, 
$$i \int_0^\infty t \ e^{ixt} \  d q(t) =i \int_0^1 t e^{ixt} \  d q(t) + i \int_1^\infty t e^{ixt} \  d q(t) $$
$$ =  f^\prime (x) +  i \int_0^1 {e^{ix /t}\over t}  \   d q (t).\eqno(2.3)$$
Recalling the relatively convergent integral from \cite{Prud}, relation (2.12.9.3) 
$$  {e^{ix /t}\over it}  =   \int_0^\infty J_0\left(2\sqrt {x y}\right)  e ^{-ity} dy,\quad x, t   > 0,\eqno(2.4)$$
where $J_0(z)$ is the Bessel function of the first kind \cite{Prud}, Vol. 2, we substitute it in (2.3). Hence after the change of the order of integration and the use of the symmetry property (1.8), we combine with (2.2) and come up with the integro-differential equation (2.1).  Our goal now is to motivate the interchange of the order of integration in the iterated integral, proving the formula
$$\int_0^1 \left( \int_0^\infty J_0\left(2\sqrt {x y}\right)  e ^{-ity} dy \right) \   d q (t) =   \int_0^\infty J_0\left(2\sqrt {x y}\right)  \left( \int_0^1 e ^{-ity} \   d q (t)\right) \  dy,\quad x >0.\eqno(2.5)$$  
To do  this, it is sufficient to justify the limit equality  
$$\lim_{Y\to \infty} \int_0^1 \left( \int_Y^\infty J_0\left(2\sqrt {x y}\right)  e ^{-ity} dy \right) \   d q (t) =0\eqno(2.6)$$
for each fixed positive $x$.   Naturally, we will appeal to the known asymptotic behavior of the Bessel function at infinity  \cite{Olver}, Section   10.17 (i)
$$J_\nu( y)=  \sqrt {{2\over \pi y}} \left[ \cos \left(y  - {\pi\nu \over 2}- {\pi\over 4} \right) -  {a (\nu) \over y}  \sin  \left( y - {\pi\nu \over 2} - {\pi\over 4} \right) + O\left({1\over y^2}\right)\right],\ y \to + \infty,\eqno(2.7)$$
where
$$a (\nu)= {\nu^2\over 2} - {1\over 8},\  \nu \in \mathbb{R}.$$
Hence, for sufficiently large $Y >0$ and $x >0,\ t \in (0,1)$, we have 
$$\int_Y^\infty J_0\left(2\sqrt {x y}\right)  e ^{-ity} dy  = {1\over \sqrt \pi x^{1/4} } \int_Y^\infty  \cos \left(2\sqrt {xy}  - {\pi\over 4} \right)  e ^{-ity} {dy\over y^{1/4}}$$
$$ +  {1\over 16\sqrt \pi x^{3/4} } \int_Y^\infty  \sin  \left(2\sqrt {xy}  - {\pi\over 4} \right)  e ^{-ity} {dy\over y^{3/4}} + O\left(Y^{-1/4}\right).\eqno(2.8)$$
As we will see from the estimates below and the finiteness of integrals (1.6) for various real $\lambda$, in order to establish the limit (2.6), it is sufficient to estimate, for instance, the integral 
$$ \int_Y^\infty  \cos \left(2\sqrt {xy}\right)\cos(ty)  {dy\over y^{1/4}},$$
because   other integrals in (2.8) can be estimated in the same manner.    With the simple substitution and  integration by parts we have 
$$\int_Y^\infty  \cos \left(2\sqrt {xy}\right) \cos(ty)  {dy\over y^{1/4}} = 2 \int_{\sqrt{Y}}^\infty  \cos \left(2 y\sqrt {x }\right) \cos(ty^2) \sqrt y \  dy$$ 
$$= - { \cos \left(2\sqrt {xY}\right) \sin\left(t Y\right) \over t Y^{1/4} }  +  {1\over 2t}  \int_{\sqrt{Y}}^\infty  \sin (ty^2) \left[ { \cos \left(2 y\sqrt {x }\right) \over  y } +  4 \sqrt x \sin(2 y\sqrt {x })  \right] \  {dy \over  \sqrt y}$$
$$= {2\sqrt x \over t}  \int_{\sqrt{Y}}^\infty  \sin (ty^2)  \sin(2 y\sqrt {x }) {dy  \over  \sqrt y} +  O\left(t^{-1} Y^{-1/4}\right).$$
Similarly,
$$ {2\sqrt x \over t} \int_{\sqrt{Y}}^\infty  \sin (ty^2)  \sin(2 y\sqrt {x }) {dy  \over  \sqrt y} = O\left(t^{-2} Y^{-3/4}\right) $$
$$+  {\sqrt x \over 2  t^2} \int_{\sqrt{Y}}^\infty  \cos (ty^2) \left[ - {3  \sin(2 y\sqrt {x }) \over  y }  +  4\sqrt{x} \cos(2 y\sqrt {x })  \right] {dy\over y^{3/2} } $$
$$=  O\left(t^{-2} Y^{-1/4}\right) .$$
Consequently, 
$$\int_0^1 \left(  \ \int_Y^\infty  \cos \left(2\sqrt {xy}\right)\cos(ty)  {dy\over y^{1/4}}   \right) d q (t) $$
$$= O\left( Y^{-1/4} \left[ \int_0^1 t^{-1} d q(t) + \int_0^1 t^{-2} d q(t)\right] \right) =   O\left( Y^{-1/4}\right),\  Y \to \infty.$$ 
Therefore, treating in the same manner other integrals from  (2.8), we  get   equality (2.6), completing the proof of Theorem 1. 
\end{proof} 

{\bf Remark 1}.   A similar to (2.1) integro-differential equation for the Fourier-Stieltjes transform (1.1)  with the derivative $f^{\prime} (x)$ inside the modified Hankel transform \cite{Tit} was exhibited in \cite{Al}. 

{\bf Corollary 1}. {\it Let $n \in \mathbb{N}.$ The values $d_n= f(2\pi n)$ and 
$$c_n= \int_0^1 t \sin(2\pi n t) dq(t)\eqno(2.9)$$
have the following integral representations in terms of the modified Hankel transform
$$ d_n=  {2\over 5}  \int_0^\infty J_0\left(2\sqrt {2\pi n y}\right) f_s(y) dy,\eqno(2.10)$$
$$ c_n=   \int_0^\infty J_0\left(2\sqrt {2\pi n y}\right) f_c(y) dy,\eqno(2.11)$$
where $f_s(x), f_c(x)$ are the Fourier-Stieltjes sine and cosine transforms  of the Minkowski question mark function $(1.10), (1.11)$,  respectively . }

\begin{proof}  Indeed, substituting in (2.1) $x=2\pi n$, we have 
$$f^\prime (2\pi n) + 2i d_n = -  \int_0^\infty J_0\left(2\sqrt {2\pi n y}\right) e^{-iy} f (y) dy.$$
In the meantime, it is not difficult to show, recalling  (1.16), that 
$$f^\prime (2\pi n)= i \int_0^1 t e^{2\pi i nt} d q(t)=  i \int_0^1 t \cos(2\pi  nt )  \ d q(t) $$
$$ - \int_0^1 t \sin (2\pi  nt )  \ d q(t) = {i\over 2} \  d_n - c_n.\eqno(2.12)$$
Hence,
$${5\over 2} \ i d_n - c_n =  - \int_0^\infty J_0\left(2\sqrt {2\pi n y}\right) e^{-iy} f (y) dy.\eqno(2.13)$$
Now taking the imaginary and real parts  of both sides of the latter equality in (2.13) with the use of (1.3), we end up with (2.10), (2.11). 
\end{proof}

Furthermore, the Salem-Zygmund theorem  \cite{Zyg}  shows that $d_n=o(1)$ implies that the Fourier-Stieltjes transform (1.1) $f(x)=o(1), |x| \to \infty$. Together with author's results in \cite{Yakeq} it leads  us to an immediate 

{\bf Corollary 2}.   {\it  Let $k \in \mathbb{N}.$  Then the Fourier-Stieltjes transforms $(1.1), (1.2)$  of the Minkowski question mark function and their consecutive derivatives  $f^{(k)}(x),  F^{(k)}(x)$  vanish at infinity.}

In particular,  Fourier-Stieltjes  coefficients $c_n$, which are defined by $(2.9)$,  vanish at infinity.   Incidentally,  it can be derived via the Fourier-Stieltjes cosine transform of the Minkowski question mark function.  In fact, owing to (1.8), (1.12)  $(x \to \infty)$, wwe find  
 
 $$ o(1)= {d\over dx} \int_0^\infty   \cos (x t )  d q(t) =  {d\over dx}  \left[ \left[ {1\over 2 - e^{ix}} + {1\over  2e^{ix}-1} \right] 
 \int_0^1   e^{ixt} d q(t) \right] $$

$$=  i  e^{ix}  \left[ {1 \over (2 - e^{ix})^2} - {2 \over  (2e^{ix}-1)^2} \right]  \int_0^1   e^{i x t }  d q(t) +  \frac{ i (e^{ix}+1) } {(2 - e^{ix}) (2e^{ix}-1)}   \int_0^1  t  e^{i x t }  d q(t).$$
Therefore,   letting $x= 2\pi n,\ n \in \mathbb{N}$, taking into account (1.16) and the value $f_s(2\pi n)=0$, we get

$$ - 2 \int_0^1  t  \sin(2\pi nt)  d q(t) = -2 c_n= o(1), \quad  n\to \infty.$$

\section{Solution to a Salem-type  problem} 

The main result is the following 

{\bf Theorem 2}. {\it   There  exists a positive bounded function $\varphi(x), \ x > 0$ such that the following linear combination of the  Fourier-Stieltjes transforms  of the Minkowski question mark function vanishes  at infinity, i.e.}  
$$ e^{-\varphi(x)} \int_0^\infty \cos\left(t\sinh\left(x+ \varphi(x)\right) \right)d q(t) +  e^{\varphi(x)} \int_0^\infty \cos\left(t\sinh\left(x- \varphi(x)\right) \right) d q(t) = o(1), \quad  x \to \infty.$$

\begin{proof}  We begin, taking (1.1) and  subtract a simple rational function and integrate by parts in the Stieltjes integral to derive

$$f (x)= \int_0^1  e^{ixt} d \left( q(t) - {2 t^2\over 1+t^2} \right) + 2 \int_0^1  e^{ixt} \  d \left( { t^2  \over 1+t^2}\right) $$

$$= - ix \int_0^1  e^{ixt} \left[  q(t) - {2 t^2\over 1+t^2} \right] dt  + O\left( {1\over x} \right),\eqno(3.1)$$
where the integrated terms vanished owing to the values $q(0)= 0,\ q(1)= 1.$ Meanwhile, passing to the limit through equality  (1.17) when $\lambda \to {\pi\over 2}-$, we get 

$${1\over \pi}\lim_{\lambda \to {\pi\over 2}-}  \int_{-\infty}^\infty \tau e^{\lambda \tau} \
\left(x+ (1+x^2)^{1/2}\right)^{i\tau} K_{i\tau}(t) d\tau  = t (1+x^2)^{1/2} \ e^{ixt}. \eqno(3.2)$$
Hence we write from (3.1)

$$f(x)+ o(1) =  {x \over  \pi i (1+x^2)^{1/2}} \int_{0}^{1}\left[
{q(t)\over t}  - {2 t\over 1+t^2} \right] $$

$$\times  \lim_{\lambda \to {\pi\over 2}-}  \int_{-\infty}^\infty \tau e^{\lambda \tau} \
\left(x+ (1+x^2)^{1/2}\right)^{i\tau} K_{i\tau}(t) d\tau\ dt,\ x \to +\infty.\eqno(3.3) $$
But since for each $x, t > 0$ and $0 \le \lambda < {\pi\over 2}$ (see (1.17) )

$$\left|\int_{-\infty}^\infty \tau e^{\lambda \tau} \ \left(x+ (1+x^2)^{1/2}\right)^{i\tau} K_{i\tau}(t) d\tau \right| \le t
\left[x+ (1+x^2)^{1/2}\right] $$
and

$$ \int_{0}^{1}\left| q(t) - {2t^2 \over 1+t^2} \right| dt \le 1 + 2 \int\limits_{0}^{1} {t^2 dt \over 1+t^2} \le  3,$$
we can take out the limit in (3.3) having the representation

$$f(x)+ o(1)  =   {x \over  \pi i (1+x^2)^{1/2}} \lim_{\lambda \to
{\pi\over 2}-} \int_{0}^{1} \left[ {q(t)\over t} - {2 t\over 1+t^2} \right] $$

$$\times  \int_{-\infty}^\infty \tau e^{\lambda \tau} \ \left(x+ (1+x^2)^{1/2}\right)^{i\tau} K_{i\tau}(t) d\tau \ dt, \ x \to +\infty.\eqno(3.4)$$
Our goal now is to invert the order of integration in (3.4). To do this we employ the uniform inequality for the modified Bessel
function (cf. \cite{LEB})
\begin{equation*}
\left|K_{i\tau}(x)\right| \le {x^{-1/4}\over \sqrt{ |\sinh (\pi\tau)|
}}, \ x > 0,\  \tau \in \mathbb{R}\backslash\{0\}
\end{equation*}
and asymptotic property  of the Minkowski question mark function near the origin. Consequently,
\begin{equation*}
\begin{split}
& \int_{0}^{1} \left|
{q(t)\over t} - {2 t\over 1+t^2}\right| \\
&\hspace{1cm}\times  \int_{-\infty}^\infty \left|\tau e^{\lambda
\tau} \ \left(x+ (1+x^2)^{1/2}\right)^{i\tau} K_{i\tau}(t)\right|
d\tau \ dt\\
&\hspace{1cm} \le \int_{0}^{1}\left|{q(t)\over t} - {2 t\over
1+t^2}\right|\ {dt\over t^{1/4}} \int_{-\infty}^\infty |\tau| {
e^{\lambda \tau}\over \sqrt{ |\sinh \pi\tau |}} d\tau < \infty, \
\lambda \in \left(0, {\pi\over 2}\right).
\end{split}
\end{equation*}
Hence by Fubini's theorem (3.4) becomes

$$ f(x)+ o(1)  =  {x \over  \pi i (1+x^2)^{1/2}} \lim_{\lambda \to
{\pi\over 2}-}  \int_{-\infty}^\infty \tau e^{\lambda \tau} \
\left(x+ (1+x^2)^{1/2}\right)^{i\tau} $$

$$ \times \int_{0}^{1} K_{i\tau}(t) \left[ {q(t)\over t} - {2 t\over 1+t^2} \right] \ dt \  d\tau, \ x \to +\infty.\eqno(3.5)$$
In the meantime, the simple change of variable and the use of the functional equation (1.5) yield 
$$  \int_0^1 K_{i\tau} (t) \left[ q(t) - {2t^2\over 1+t^2} \right] {dt \over t} =  
 \int_1^\infty  K_{i\tau} \left({1\over t} \right) \left[ q \left({1\over t} \right) - {2\over 1+t^2} \right] {dt\over t}$$

$$=  -  \int_1^\infty  K_{i\tau} \left({1\over t} \right) \left[ q(t) - {2t^2 \over 1+t^2} \right] {dt\over t}.$$
Hence in the same manner 

$$  \int_0^1 \left[ K_{i\tau} (t)  -   K_{i\tau} \left({1\over t} \right) \right] \left[ q(t) - {2 t^2\over 1+t^2} \right] {dt\over t}$$

$$ =  -  \int_0^\infty  K_{i\tau} \left({1\over t} \right) \left[ q(t) - {2 t^2\over 1+t^2} \right] {dt\over t} $$

$$=  \int_0^\infty  K_{i\tau} (t) \left[ q(t) - {2 t^2\over 1+t^2} \right] {dt\over t} .\eqno(3.6)$$

But via  (3.2) and (1.5)

$$-  {x\over \pi i (1+x^2)^{1/2}} \lim_{\lambda \to \pi/2 -}\int_{-\infty}^\infty \tau e^{\lambda \tau} \left(x+ (1+x^2)^{1/2}\right)^{i\tau} \int_0^1 K_{i\tau} \left({1\over t} \right)  \left[ q(t) - {2t^2\over 1+t^2} \right] {dt d\tau\over t} $$

$$= ix \int_0^1 e^{i x/t}  \left[ q(t) - {2t^2\over 1+t^2} \right] {dt \over t^2} =  ix \int_1^\infty e^{it x}  \left[ q \left({1\over t} \right) - {2 \over 1+t^2} \right]  dt$$

$$=  - ix \int_1^\infty e^{it x}  \left[ q(t) - {2 t^2 \over 1+t^2} \right]  dt =  \int_1^\infty e^{it x}  d q(t) +  o(1), \ x \to + \infty.\eqno(3.7)$$
However  from (1.12) we find

$$ \int_1^\infty e^{it x}  d q(t) = {e^{ix}\over 2 - e^{ix}}  \int_0^1  e^{it x}  d q(t) .$$
Therefore, combining with (3.5), (3.6), (3.7), we deduce the equality 

$$ {2\ f(x) \over 2 - e^{ix}}  + o(1) =  {x\over \pi i (1+x^2)^{1/2}} \lim_{\lambda \to \pi/2 -}\int_{-\infty}^\infty \tau e^{\lambda \tau} (x+ (x^2+1)^{1/2})^{i\tau} $$

$$\times \int_0^\infty K_{i\tau} (t) \left[ q(t) - {2 t^2\over 1+t^2} \right] {dt d\tau \over t} , \ x \to +\infty.\eqno(3.8)$$
Meanwhile,  from (1.5) and asymptotic behavior of the Minkowski question mark function we get

$$ \left[ q(t) - {2 t^2\over 1+t^2}\right] {1\over t} =  O( t),\ t \to 0,$$ 

$$ \left[ q(t) - {2 t^2\over 1+t^2}\right] {1\over t}  = O\left( { 2\over t}  - {2 t \over 1+t^2} \right) = O( t^{-3} ), \ t \to +\infty,$$
and therefore

$$ \left[ q(t) - {2 t^2\over 1+t^2}\right] {1\over t} \in L_p\left(\mathbb{R}_+\right),\quad p \ge 1.$$
This circumstance together with the asymptotic behavior of the modified Bessel function  $K_{i\tau} (t)$  at zero and infinity \cite{Olver} for fixed $\tau$ permits us to apply the Parseval equality for the  Fourier  cosine transform (see \cite{Tit}, Theorem 52). Thus employing (1.18), we obtain 
$$ \int_0^\infty K_{i\tau} (t) \left[ q(t) - {2 t^2\over 1+t^2} \right] {dt\over t} = {1\over  \cosh( \pi \tau/2) } \int_0^\infty   {\cos( \tau\ \log(y+(1+y^2)^{1/2} )) \over (y^2+1)^{1/2}} $$

$$\times  \int_0^\infty \cos(t y) \left[ q(t) - {2t^2\over 1+t^2} \right] {dt\ dy\over t}.\eqno(3.9)$$
Substituting the right-hand side of (3.9) into the right-hand side of (3.8) and making simple substitutions, we derive 

$$ {x\over \pi i (1+x^2)^{1/2}} \int_{-\infty}^\infty \tau e^{\lambda \tau} (x+ (x^2+1)^{1/2})^{i\tau}  \int_0^\infty K_{i\tau} (t) \left[ q(t) - {2 t^2\over 1+t^2} \right] {dt d\tau \over t}$$

$$=  {x\over \pi i (1+x^2)^{1/2}}  \int_{-\infty}^\infty  {\tau \ e^{\lambda \tau} \over  \cosh( \pi \tau/2) } \  (x+ (x^2+1)^{1/2})^{i\tau} $$

$$\times  \int_0^\infty  \cos( \tau y)  \int_0^\infty \cos(t \sinh(y) ) \left[ q(t) - {2t^2\over 1+t^2} \right] {dt\ dy\over t} \ d\tau.\eqno(3.10)$$
Moreover,  the $L_2$-theory of the Fourier transform says (see \cite{Tit} )

$$G(y)= \int_0^\infty \cos(t y ) \left[ q(t) - {2t^2\over 1+t^2} \right] {dt\  \over t}  \in L_2\left(\mathbb{R}_+\right)$$
and therefore $G(\sinh(y)) \in L_2\left(\mathbb{R}_+\right).$ This fact yields

$$ \int_0^\infty  \cos( \tau y) \ G(\sinh(y)) dy  \in L_2\left(\mathbb{R}\right)$$
as a function of $\tau$. Therefore applying the Cauchy- Schwarz inequality, we easily establish the absolute and uniform convergence with respect to $x$ of the integral by $\tau$ on the right-hand side of (3.10) for each $\lambda \in [ 0, \pi/2)$.  Hence it is possible to differentiate under the integral sign on the right-hand side of (3.10),  and we obtain the equality 

$${x\over \pi i (1+x^2)^{1/2}}  \int_{-\infty}^\infty  {\tau \ e^{\lambda \tau} \over  \cosh( \pi \tau/2) } \  (x+ (x^2+1)^{1/2})^{i\tau} $$

$$\times  \int_0^\infty  \cos( \tau y)  \int_0^\infty \cos(t \sinh(y) ) \left[ q(t) - {2t^2\over 1+t^2} \right] {dt\ dy\over t}\ d\tau$$

$$= -  {x\over \pi} \  {d\over dx}  \int_{-\infty}^\infty  { e^{\lambda \tau} \over  \cosh( \pi \tau/2) } \  (x+ (x^2+1)^{1/2})^{i\tau} $$

$$\times  \int_0^\infty  \cos( \tau y)  \int_0^\infty \cos(t \sinh(y) ) \left[ q(t) - {2t^2\over 1+t^2} \right] {dt\ dy\over t}\ d\tau.\eqno(3.11)$$
Moreover, the latter iterated integral can be treated via the generalized Parseval equality for Fourier transform (see \cite{Tit}, Theorem 64). An alternative approach is to use Fubini's theorem and appeal to formula (1.19). Consequently, employing (1.19), it becomes 

$$  \int_{-\infty}^\infty  { e^{\lambda \tau} \over  \cosh( \pi \tau/2) } \  (x+ (x^2+1)^{1/2})^{i\tau} $$

$$\times  \int_0^\infty  \cos( \tau y)  \int_0^\infty \cos(t \sinh(y) ) \left[ q(t) - {2t^2\over 1+t^2} \right] {dt\ dy\over t}\ d\tau$$

$$= \int_{-\infty}^\infty  \left[  \cos ( \lambda + i( y+ \log (x+ (x^2+1)^{1/2}) ))\right]^{-1}$$

$$\times  \int_0^\infty \cos(t \sinh (y) ) \left[ q(t) - {2 t^2\over 1+t^2} \right] {dt\ dy\over t}.\eqno(3.12)$$
Then, fixing a  positive $\delta$, we split the latter integral with respect to $y$ as follows

$$ \int_{-\infty}^\infty  \left[  \cos ( \lambda + i( y+ \log (x+ (x^2+1)^{1/2}) ))\right]^{-1} \int_0^\infty \cos(t \sinh (y) ) \left[ q(t) - {2t^2\over 1+t^2} \right] {dt\ dy\over t}$$

$$= \left(  \int_{-\infty }^{-\log (x+ (x^2+1)^{1/2})-\delta} +  \int_{- \log (x+ (x^2+1)^{1/2})-\delta}^ {- \log (x+ (x^2+1)^{1/2})+\delta} + \int_{- \log (x+ (x^2+1)^{1/2})+\delta}^\infty \right) $$

$$\times \left[  \cos ( \lambda + i( y+ \log (x+ (x^2+1)^{1/2}) ))\right]^{-1}  \int_0^\infty \cos(t \sinh (y) ) \left[ q(t) - {2t^2\over 1+t^2} \right] {dt\ dy\over t}.\eqno(3.13)$$
Considering the first integral on the right-hand side of the latter equality, we can differentiate it under the integral sign with respect to $x$ and then pass to the limit when $\lambda \to \pi/2 -$ owing to the absolute and uniform convergence. Thus it gives 

$$-  {x\over \pi  } \lim_{\lambda \to \pi/2 -} {d\over dx}  \int_{-\infty }^{-\log (x+ (x^2+1)^{1/2})-\delta} \left[  \cos ( \lambda + i( y+ \log (x+ (x^2+1)^{1/2}) ))\right]^{-1} $$

$$\times  \int_0^\infty \cos(t \sinh (y) ) \left[ q(t) - {2t^2\over 1+t^2} \right] {dt\ dy\over t}$$

$$=  {x\over \pi i   (x^2+1)^{1/2} \sinh(\delta)}  \int_0^\infty \cos(t \sinh ( \log (x+ (x^2+1)^{1/2})+\delta) ) \left[ q(t) - {2t^2\over 1+t^2} \right] {dt \over t}$$

$$- {x\over \pi  i  (x^2+1)^{1/2} }   \int_{-\infty}^{-\log (x+ (x^2+1)^{1/2})-\delta}  { \cosh (  y+ \log (x+ (x^2+1)^{1/2}) ) \over \sinh^2 ( y+ \log (x+ (x^2+1)^{1/2}) ) } $$

$$\times  \int_0^\infty \cos(t \sinh (y) ) \left[ q(t) - {2t^2\over 1+t^2} \right] {dt\ dy\over t}.\eqno(3.14)$$
However, since
$$\int_0^\infty \left|  q(t) - {2t^2\over 1+t^2} \right| {dt \over t} < \infty$$
Therefore owing to the Riemann-Lebesgue lemma 

$${x\over \pi i   (x^2+1)^{1/2} \sinh(\delta)}  \int_0^\infty \cos(t \sinh ( \log (x+ (x^2+1)^{1/2})+\delta) ) \left[ q(t) - {2t^2\over 1+t^2} \right] {dt \over t} = o(1),\ x \to + \infty.$$
Moreover, making a simple substitution we arrive at the estimate

$$ \left| \int_{-\infty}^{-\log (x+ (x^2+1)^{1/2})-\delta}  { \cosh (  y+ \log (x+ (x^2+1)^{1/2}) ) \over \sinh^2 ( y+ \log (x+ (x^2+1)^{1/2}) ) } \right.$$

$$\times \left.  \int_0^\infty \cos(t \sinh (y) ) \left[ q(t) - {2t^2\over 1+t^2} \right] {dt\ dy\over t}\right| $$

$$= \left| \int_{ \delta}^\infty   { \cosh ( y ) \over \sinh^2 ( y) }    \int_0^\infty \cos(t \sinh (y +\log (x+ (x^2+1)^{1/2} ) ) \left[ q(t) - {2t^2\over 1+t^2} \right] {dt \ dy \over t}\right| $$

$$ \le  \int_{ \delta}^\infty   { \cosh ( y) \ dy \over \sinh^2 ( y) }  \int_0^\infty \left|  q(t) - {2t^2\over 1+t^2}  \right| {dt \over t}  < \infty,$$
which allows to pass to the limit under the integral sign in the integral with respect to $y$ via the dominated convergence theorem. Thus recalling the Riemann - Lebesgue lemma, we establish the following asymptotic equality for the first integral  on the right-hand side of (3.13), namely

$$ -  {x\over \pi  } \lim_{\lambda \to \pi/2 -} {d\over dx}  \int_{-\infty }^{-\log (x+ (x^2+1)^{1/2})-\delta} \left[  \cos ( \lambda + i( y+ \log (x+ (x^2+1)^{1/2}) ))\right]^{-1} $$

$$\times  \int_0^\infty \cos(t \sinh (y) ) \left[ q(t) - {2t^2\over 1+t^2} \right] {dt\ dy\over t} = o(1),\ x \to +\infty.$$
Analogously,  we treat the third integral on the right-hand side of (3.13) to derive   

$$ -  {x\over \pi  } \lim_{\lambda \to \pi/2 -} {d\over dx}  \int_{-\log (x+ (x^2+1)^{1/2})+\delta}^\infty  \left[  \cos ( \lambda + i( y+ \log (x+ (x^2+1)^{1/2}) ))\right]^{-1} $$

$$\times  \int_0^\infty \cos(t \sinh (y) ) \left[ q(t) - {2t^2\over 1+t^2} \right] {dt\ dy\over t} = o(1),\ x \to +\infty.$$
Finally, let us consider the middle integral  on the right-hand side of (3.13).  We write
$$-  {x\over \pi  }  {d\over dx}  \int_{-\log (x+ (x^2+1)^{1/2})-\delta}^{-\log (x+ (x^2+1)^{1/2})+\delta} \left[  \cos ( \lambda + i( y+ \log (x+ (x^2+1)^{1/2}) ))\right]^{-1} $$

$$\times  \int_0^\infty \cos(t\sinh (y) ) \left[ q(t) - {2t^2\over 1+t^2} \right] {dt\ dy\over t}$$

$$= -  {x\over \pi  }  {d\over dx}  \int_{-\delta}^{\delta} \left[  \cos ( \lambda + i y)\right]^{-1} $$

$$\times  \int_0^\infty \cos(t \sinh (y- \log (x+ (x^2+1)^{1/2})) ) \left[ q(t) - {2 t^2\over 1+t^2} \right] {dt\ dy\over t}$$

$$= -   {x\over \pi  (x^2+1)^{1/2} } \int_{-\delta}^{\delta}  {\cosh (y-\log (x+ (x^2+1)^{1/2})) \over \cos ( \lambda + i y) } $$

$$\times   \int_0^\infty \sin (t \sinh (y-\log (x+ (x^2+1)^{1/2})) ) \left[ q(t) - {2t^2\over 1+t^2} \right] dt dy,\eqno(3.15)$$
where the differentiation  under the integral sign in (3.15) is permitted  because $\lambda \in [0, \pi/2)$ and

$$ q(t) - {2t^2\over 1+t^2} \in L_1(\mathbb{R}_+).$$
Further,

$$ -   {x\over \pi  (x^2+1)^{1/2} } \int_{-\delta}^{\delta}  {\cosh (y-\log (x+ (x^2+1)^{1/2})) \over \cos ( \lambda + i y) } $$

$$\times   \int_0^\infty \sin (t \sinh (y-\log (x+ (x^2+1)^{1/2})) ) \left[ q(t) - {2t^2\over 1+t^2} \right] dt dy$$

$$= -   {x\over \pi  } \int_{-\delta}^{\delta}  {\cosh (y)  \over \cos ( \lambda + i y) } $$

$$\times   \int_0^\infty \sin (t \sinh (y-\log (x+ (x^2+1)^{1/2})) ) \left[ q(t) - {2t^2\over 1+t^2} \right] dt dy $$

$$+  {x^2\over \pi (x^2+1)^{1/2} } \int_{-\delta}^{\delta}  {\sinh (y)  \over \cos ( \lambda + i y) } $$

$$\times   \int_0^\infty \sin (t \sinh (y-\log (x+ (x^2+1)^{1/2})) ) \left[ q(t) - {2t^2\over 1+t^2} \right] dt dy .$$
However, the latter integral converges uniformly in $\lambda \in [ \pi/2 - \varepsilon, \ \pi/2],$ where $\ \varepsilon > 0$ is a small fixed number. Indeed, we have the estimate

$$ \int_{-\delta}^{\delta} \left| {\sinh (y)  \over  \cos ( \lambda + i y) } \right| $$

$$\times  \left| \int_0^\infty \sin (t \sinh (y-\log (x+ (x^2+1)^{1/2})) ) \left[ q(t) - {2t^2\over 1+t^2} \right] dt dy \right| $$

$$ \le 2\delta \csc (\lambda)  \int_0^\infty \left| q(t) - {2t^2\over 1+t^2} \right| dt < \infty. $$
Therefore passing to the limit under the integral sign when $\lambda \to \pi/2-$ and then integrating by parts, we obtain 

$${x^2\over \pi (x^2+1)^{1/2} }  \lim_{\lambda \to \pi/2 -} \int_{-\delta}^{\delta}  {\sinh (y)  \over \cos ( \lambda + i y) } $$

$$\times   \int_0^\infty \sin (t \sinh (y-\log (x+ (x^2+1)^{1/2})) ) \left[ q(t) - {2t^2\over 1+t^2} \right] dt dy $$

$$= - {x^2\over \pi  i (x^2+1)^{1/2} }   \int_{-\delta}^{\delta}   \int_0^\infty \sin (t \sinh (y-\log (x+ (x^2+1)^{1/2})) ) \left[ q(t) - {2t^2\over 1+t^2} \right] dt dy $$

$$= {x^2\over \pi i  (x^2+1)^{1/2}} [ \cosh (\delta -\log (x+ (x^2+1)^{1/2})) ]^{-1} $$

$$\times \int_0^\infty \cos (t \sinh (\delta -\log (x+ (x^2+1)^{1/2})) ) \left[ q(t) - {2t^2\over 1+t^2} \right] {dt\over t} $$

$$ - {x^2\over \pi i  (x^2+1)^{1/2}} [ \cosh (\delta +\log (x+ (x^2+1)^{1/2}))]^{-1}$$

$$\times  \int_0^\infty \cos (t \sinh ( \delta+\log (x+ (x^2+1)^{1/2})) ) \left[ q(t) - {2t^2\over 1+t^2} \right] {dt\over t} $$

$$ + {x^2\over \pi i  (x^2+1)^{1/2}}  \int_{-\delta}^{\delta} { \tanh (y -\log (x+ (x^2+1)^{1/2}))\over  \cosh (y -\log (x+ (x^2+1)^{1/2})) }$$

$$\times \int_0^\infty \cos (t \sinh (y- \log (x+ (x^2+1)^{1/2})) ) \left[ q(t) - {2t^2\over 1+t^2} \right] {dt\over t} \ dy = o(1),\  x \to + \infty$$
due to the dominated convergence theorem and the Riemann- Lebesgue lemma.   Further,  it is easily seen that

$$- {x\over \pi} \ \int_{-\delta}^{\delta}  {\cosh (y)  \over \cos ( \lambda + i y) } $$

$$\times   \int_0^\infty \sin (t \sinh (y-\log (x+ (x^2+1)^{1/2})) ) \left[ q(t) - {2t^2\over 1+t^2} \right] dt dy $$

$$=   {2 x\over \pi}   \int_{0}^{\delta} { \cos (\lambda) \cosh^2(y) \over \sinh^2(y)+ \cos^2(\lambda) } 
 \int_0^\infty \sin(x t \cosh(y)) \cos(t  (x^2+1)^{1/2} \sinh (y) )$$
 
 $$\times   \left[ q(t) - {2t^2\over 1+t^2} \right] dt dy $$

$$-    { 2i x\over \pi} \sin (\lambda) \int_{0}^{\delta} { \sinh(y) \cosh(y) \over \sinh^2(y)+ \cos^2(\lambda) } 
 \int_0^\infty \cos (x t \cosh(y)) \sin (t  (x^2+1)^{1/2} \sinh (y) )$$
 
 $$\times   \left[ q(t) - {2t^2\over 1+t^2} \right] dt\ dy,\  \lambda \in \left[0, {\pi\over 2}\right). \eqno(3.16)$$
Hence we split  the first integral on the right-hand side of the latter equality in the following manner   

$$  {2 x\over \pi}  \int_{0}^{\delta} {  \cos (\lambda) \cosh^2(y) \over \sinh^2(y)+ \cos^2(\lambda) } 
 \int_0^\infty \sin(x t \cosh(y)) \cos(t  (x^2+1)^{1/2} \sinh (y) )$$
 
 $$\times   \left[ q(t) - {2t^2\over 1+t^2} \right] dt dy $$

$$=   { 2 x\over \pi}  \int_{0}^{\delta } { \cos (\lambda) \sinh^2(y) \over \sinh^2(y)+ \cos^2(\lambda) } 
 \int_0^\infty \sin(x t \cosh(y)) \cos(t  (x^2+1)^{1/2} \sinh (y) )$$
 
 $$\times   \left[ q(t) - {2t^2\over 1+t^2} \right] dt dy $$

$$ +{ 2x\over \pi}  \int_{0}^{ \delta/ \cos (\lambda)} {\cos^2 (\lambda)  \over \sinh^2(y\cos(\lambda))+ \cos^2(\lambda) } 
 \int_0^\infty \sin(x t \cosh(y\cos(\lambda))) \cos(t  (x^2+1)^{1/2} \sinh (y\cos(\lambda)) )$$
 
 $$\times   \left[ q(t) - {2t^2\over 1+t^2} \right] dt dy .$$
But, evidently,  via the dominated convergence theorem

$$  {2 x\over \pi}   \lim_{\lambda \to \pi/2 -}\ \int_{0}^{\delta} { \cos (\lambda)  \sinh^2(y) \over \sinh^2(y)+ \cos^2(\lambda) } 
 \int_0^\infty \sin(x t \cosh(y)) \cos(t  (x^2+1)^{1/2} \sinh (y) )$$
 
 $$\times   \left[ q(t) - {2t^2\over 1+t^2} \right] dt dy = 0.$$
Moreover, since

$$\int_{0}^{ \delta/ \cos (\lambda)} {\cos^2 (\lambda)  \over \sinh^2(y\cos(\lambda))+ \cos^2(\lambda) } 
 \int_0^\infty \left| \sin(x t \cosh(y\cos(\lambda))) \cos(t  (x^2+1)^{1/2} \sinh (y\cos(\lambda)) )\right|$$
 
 $$\times   \left| q(t) - {2t^2\over 1+t^2} \right| dt dy $$

$$\le  \int_{0}^{ \infty}  { dy \over y^2+1 }  \int_0^\infty  \left| q(t) - {2t^2\over 1+t^2} \right| dt < \infty,$$
we get,  owing to the same arguments 

$$  {2 x\over \pi}   \lim_{\lambda \to \pi/2 -}  \int_{0}^{ \delta} {\cos (\lambda)  \over \sinh^2(y)+ \cos^2(\lambda) }  \int_0^\infty \sin(x t \cosh(y)) \cos(t  (x^2+1)^{1/2} \sinh (y) )$$
 
 $$\times   \left[ q(t) - {2t^2\over 1+t^2} \right] dt dy $$

$$=  x \int_0^\infty \sin(x t )  \left[ q(t) - {2t^2\over 1+t^2} \right] dt =   \int_0^\infty \cos (x t ) \ d  \left( q(t) - {2t^2\over 1+t^2} \right). \eqno(3.17)$$
Now,   returning to (3.16), we treat the latter integral on its right-side, employing  equality (1.20).  Hence we deduce

$$ - {2 i x\over \pi}  \sin (\lambda) \int_{0}^{\delta} { \sinh(y) \cosh(y) \over \sinh^2(y)+ \cos^2(\lambda) } 
 \int_0^\infty \cos (x t \cosh(y)) \sin (t  (x^2+1)^{1/2} \sinh (y) )$$
 
 $$\times   \left[ q(t) - {2t^2\over 1+t^2} \right] dt\ dy$$

$$=   -  i x   \sin (\lambda)  \int_0^\infty  e^{ - t  (x^2+1)^{1/2}  \cos (\lambda) } \cos (x t \sin (\lambda) ) \left[ q(t) - {2t^2\over 1+t^2} \right] dt$$
$$+  {2 i x\over \pi}   \sin (\lambda)  \int_0^\infty \int_{\sinh(\delta) }^\infty  { u    \cos (x t (u^2+1)^{1/2} ) \sin (t  (x^2+1)^{1/2}  u ) \over  u^2 + \cos^2(\lambda) }  $$
 
 $$\times   \left[ q(t) - {2t^2\over 1+t^2} \right]  du dt.\eqno(3.18)$$
Then

$$ -  i x    \lim_{\lambda \to \pi/2 -}\  \sin (\lambda)  \int_0^\infty  e^{ - t  (x^2+1)^{1/2}  \cos (\lambda) } \cos (x t \sin (\lambda) ) \left[ q(t) - {2t^2\over 1+t^2} \right] dt$$

$$= -  i x    \int_0^\infty   \cos (x t ) \left[ q(t) - {2t^3\over 1+t^3} \right] dt =   i   \int_0^\infty   \sin  (x t )\  d \left( q(t) - {2t^2\over 1+t^2} \right) .\eqno(3.19)$$
Therefore, returning to the previous results, we derive the following asymptotic equality from (3.8), (3.17), (3.18), (3.19)  

$$  {2 i x\over \pi}    \lim_{\lambda \to \pi/2 -}  \sin (\lambda)  \int_0^\infty \int_{\sinh(\delta) }^\infty  { u    \cos (x t (u^2+1)^{1/2} ) \sin (t  (x^2+1)^{1/2}  u ) \over  u^2 + \cos^2(\lambda) }  $$
 
 $$\times   \left[ q(t) - {2t^2\over 1+t^2} \right]  du dt = o (1),\quad x \to +\infty.\eqno(3.20)$$
Hence, taking two different $0 < \delta_1 < \delta_2$,  equality  (3.20) yields when $x \to \infty$

$$    {2 i x\over \pi}    \lim_{\lambda \to \pi/2 -}   \int_0^\infty \int_{\sinh(\delta_1) }^ {\sinh(\delta_2)}   { u    \cos (x t (u^2+1)^{1/2} ) \sin (t  (x^2+1)^{1/2}  u ) \over  u^2 + \cos^2(\lambda) }    \left[ q(t) - {2t^2\over 1+t^2} \right]  du dt =  o(1).$$
 Passing to the limit under the integral sign on the left-hand side of the latter equality via the dominated convergence theorem, we change the order of integration by Fubini's theorem and  make  simple substitutions to  write  it in the form
 
$$    \sinh(x)   \int_{\delta_1 }^ {\delta_2}  \coth(u)  \int_0^\infty  \left[   \sin (t \sinh(x +u)  ) +  \sin (t \sinh(x - u) \right]  \left[ q(t) - {2t^2\over 1+t^2} \right] dt du =  o(1).\eqno(3.21)$$
The integration by parts in the integral with respect to $t$ in (3.21) and elimination of small terms imply

 $$    \sinh(x)   \int_{\delta_1 }^ {\delta_2}  {\coth(u) \over \sinh(x +u) }  \left(\int_0^\infty    \cos (t \sinh(x +u)  ) \  d q(t) \right) du$$
 
 $$+  \sinh(x)   \int_{\delta_1 }^ {\delta_2}  {\coth(u) \over \sinh(x -u) } \left( \int_0^\infty    \cos (t \sinh(x -u)  )\   d q(t)\right) du =  o(1),\ x \to \infty.\eqno(3.22)$$
 Since the Fourier-Stieltjes integrals in (3.22) are continuous functions of $u$, we employ the mean value  theorem for the corresponding integral to get 

 $$    \sinh(x) \coth(\varphi(x)) (\delta_2-\delta_1) \left[  {1 \over \sinh(x + \varphi(x)) }  \int_0^\infty    \cos (t \sinh(x + \varphi(x))  ) \  d q(t) \right.$$
 
 $$\left. +   {1 \over \sinh(x - \varphi(x)) } \  \int_0^\infty    \cos (t \sinh(x - \varphi (x))  )\   d q(t) \right] =  o(1),\ x \to \infty,$$
 where $\delta_1 \le \varphi(x) \le \delta_2.$  Thus, eliminating again small terms, we obtain

 $$   e^{-\varphi(x)}  \int_0^\infty    \cos (t \sinh(x + \varphi(x))  ) \  d q(t) $$
 
 $$ +   e^{\varphi(x)} \int_0^\infty    \cos (t \sinh(x - \varphi (x))  )\   d q(t) =  o(1),\ x \to \infty,$$
completing the proof of Theorem 2.

\end{proof}

Finally, denoting by
$$d_{n,m} = \int_0^1 \cos(2\pi n t) \ d q^m(t),\ m \in \mathbb{N}\eqno(3.23)$$
the Fourier-Stieltjes coefficients of the power $q^m(t)$, we establish 

{\bf Corollary 3}.  {\it  Let $m \in \mathbb{N}$. The following asymptotic relations hold valid

$$d_{n,2m} -  {1\over 2}  \sum_{k=2}^{2m-1} (-1)^{k+1} \binom {2m} {k}   d_{n,k} =  o(1), \quad  n\to \infty,\eqno(3.24)$$

$$\sum_{k=2}^{2(m-1)} (-1)^{k+1} \binom {2m-1} {k}   d_{n,k} =  o(1), \quad  n\to \infty.\eqno(3.25)$$
In particular,  $d_{n,2}= o(1),\ n \to \infty$}.

\begin{proof}   Indeed,  recalling  (1.3), we  have the equality 
$$d_{n,2m} = \sum_{k=1}^{2m} (-1)^k \binom {2m} {k}    \int_0^1 \cos(2\pi n t)\  d  q^k(1-t) = \sum_{k=1}^{2m} (-1)^{k+1} \binom {2m} {k}   d_{n,k}.$$ 
Hence,

$$d_{n,2m} = {1\over 2}  \sum_{k=1}^{2m-1} (-1)^{k+1} \binom {2m} {k}   d_{n,k},\ m \in \mathbb{N}.$$
But  $d_{n,1} \equiv d_n= o(1),\ n \to \infty$. This gives (3.24).  Analogously,

$$ d_{n,2m-1} = \sum_{k=1}^{2m-1} (-1)^k \binom {2m-1} {k}    \int_0^1 \cos(2\pi n t)\  d  q^k(1-t) = \sum_{k=1}^{2m-1} (-1)^{k+1} \binom {2m-1} {k}   d_{n,k},$$ 
i.e.

$$  \sum_{k=1}^{2(m-1)} (-1)^{k+1} \binom {2m-1} {k}   d_{n,k} =0.$$ 
The latter equality yields (3.25).  

\end{proof}

{\bf Open problem}. { Prove or disprove}

$$ d_{n,m}= o(1),\quad  n \to \infty,\quad m \ge 3.$$

\noindent {{\bf Acknowledgments}}\\

The  work  was partially supported by CMUP (UID/MAT/00144/2013), which is funded by FCT (Portugal) with national (MEC) and European structural funds through the programs FEDER, under the partnership agreement PT2020.   The author is sincerely indebted to the referee for a careful reading of the manuscript and useful comments, which rather improved the presentation of the paper.  \\

\vspace{1cm}

\noindent {\sc Semyon Yakubovich}, Department of Mathematics, Faculty of Sciences, University of Porto,
Campo Alegre st., 687,  4169-007 Porto, Portugal. {\tt syakubov@fc.up.pt}\\

\end{document}